\documentclass[graybox]{svmult}

\usepackage{amsmath}
\usepackage[english]{babel}
\usepackage{txfonts}
\usepackage{helvet}
\usepackage{soul}   
\usepackage{courier}        
\usepackage{type1cm}        
\usepackage{makeidx}         
\usepackage{graphicx}        
\usepackage[bottom]{footmisc}

\usepackage{subfigure}
\usepackage[abs]{overpic}
\makeindex             

\def\NN{\hbox{I\kern-.2em\hbox{N}}}
\def\RR{{\mathop{{\rm I}\kern-.2em{\rm R}}\nolimits}}

\def\Q{{\bf Q}}

\def\d{{\bf d}}
\def\F{{\bf F}}
\def\g{{\bf g}}
\def\p{{\bf p}}
\def\q{{\bf q}}

\def\s{{\bf s}}
\def\t{{\bf t}}
\def\x{{\bf x}}
\def\y{{\bf y}}

\def\bfax{\mbox{\boldmath$\alpha$}}

\def\bnu{\mbox{\boldmath$\nu$}}

\definecolor{verde}{RGB}{46,139,87}

\newcommand{\be}{\begin{equation}}
\newcommand{\ee}{\end{equation}}
\newcommand{\ba}{\begin{eqnarray}}
\newcommand{\ea}{\end{eqnarray}}

\begin{document}

\title*{A collocation IGA-BEM for 3D potential problems on unbounded domains}
\titlerunning{3D IGA-BEM}

\author{Antonella Falini, Carlotta Giannelli, Tadej Kandu\v{c}, M. Lucia Sampoli and Alessandra Sestini}
\institute{Antonella Falini, Dipartimento di Informatica, Universit\`a degli Studi di Bari Aldo Moro, Italy \email{antonella.falini@uniba.it} \\
\and Carlotta Giannelli, Dipartimento di Matematica e Informatica Ulisse Dini, Universit\`a degli Studi di Firenze, Italy \email{carlotta.giannelli@unifi.it}\\
\and  Tadej Kandu\v{c}, Faculty of Mathematics and Physics, University of Ljubljana, Ljubljana, Slovenia \email{tadej.kanduc@fmf.uni-lj.si}\\
\and M. Lucia Sampoli, Dipartimento di Ingegneria dell'Informazione e Scienze Matematiche, Universit\`a degli Studi di Siena, Italy\email{marialucia.sampoli@unisi.it}\\
\and Alessandra Sestini, Dipartimento di Matematica e Informatica Ulisse Dini, Universit\`a degli Studi di Firenze, Italy \email{alessandra.sestini@unifi.it}}
\authorrunning{Falini, Giannelli, Kandu\v{c}, Sampoli, Sestini}
%
%
\maketitle

\abstract{In this paper the numerical solution of potential problems defined on 3D unbounded domains is addressed with Boundary Element Methods (BEMs), since in this way the problem is studied only on the boundary, and thus any finite approximation of the infinite domain can be avoided. The isogeometric analysis (IGA) setting is considered and in particular B-splines and NURBS functions are taken into account. In order to exploit all the possible benefits from using spline spaces, an important point is the development of specific cubature formulas for weakly and nearly singular integrals. Our proposal for this aim is based on spline quasi-interpolation and on the use of a spline product formula. Besides that, a robust singularity extraction procedure is introduced as a preliminary step  and an efficient function-by-function assembly phase  is adopted.  A selection of numerical examples confirms that the numerical solutions reach the expected convergence orders.
}


\section{Introduction}\label{sec:intro}

The framework of Isogeometric Analysis (IGA), introduced in 2005 by Thomas J.~R.~Hughes and coauthors, has emerged as a very attractive approach to solve differential problems \cite{LibroHughes, Hughes_2005}. IGA bridges the
gap between Computational Mechanics and Computer Aided Design (CAD) by employing a unified representation of geometric objects and  physical quantities.  In particular, the spline basis functions commonly used for CAD geometry parameterization are also used to span the discretization space adopted for the numerical solution. This not only eliminates the need for slow and error-prone conversion processes between the geometry and simulation models, but  also enables the possibility of a widespread use of design optimization tools in simulation-based engineering.

While IGA has achieved great success within the finite element paradigm by attracting most of the IGA researchers, a remaining
critical challenge is how to efficiently and exactly parameterize a volumetric domain from its boundary \cite{al2016Brep,cohen_etal2010,pan2020volumetric}. Three dimensional problems need trivariate NURBS solids for analysis, but CAD systems use a boundary representation where only bivariate NURBS are involved. A possible solution to this issue consists in adopting boundary element methods (BEMs) \cite{beer2020LIBRO, BEMbook}, whose isogeometric formulation is briefly referred to as IGA-BEM and has recently been succesfully applied to several classical problems, e.g., Laplace and Helmholtz equation \cite{TauRodHug, Dolz18}, elastostatics and elasticity \cite{Simp2012, Wang15, beer_etal2017, nguyen2017}, potential flow \cite{gong2017} and wave-resistance problems \cite{ginnis}.
 BEMs reformulate the original differential problem through the fundamental solution  into boundary integral equations. Thus, adopting a boundary representation,  they can be easily coupled with standard CAD techniques  and ensure dimension reduction of the computational domain where the numerical simulation is developed. Besides these key points, they are attractive because offering an easy treatment  of problems on unbounded domains, see for example the potential or acoustic problems considered in \cite{Coox, Nguyen16}. Within this kind of problems, screen or crack problems are particularly important for applications, being in this case the domain given by the whole Euclidean space external to an open limited arc in 2D or surface in 3D.

Relying on the general NURBS parametric description of the domain boundary and on (refinements of) the related spline space for the approximation of the unknowns, in this paper we confirm  the effectiveness of the IGA-BEM in the 3D setting, previously experimented  in 2D \cite{ADSS1}.  In particular, we focus on a collocation discretization strategy for 3D exterior potential problems. The combination of the superior smoothness of spline basis functions with the low computational cost and simplicity  of  collocation  techniques  seems  to  constitute  an optimal  basis  for  accurately  modeling complex and computationally demanding  problems, at least when smooth geometries are dealt with,  as also outlined in  \cite{TauRodHug}.  However it is worth mentioning that considering a Galerkin based variant of the developed method   could be of great interest in the future, either because the method could become more robust for treating problems with low regular solutions \cite{sutradhar2008SGBEM} or even from a more theoretical point of view, since most of the convergence results available in the literature refer to the Galerkin approach.\\
Since an effective implementation of BEMs surely requires suitable rules for the approximation of arising singular integrals, in this paper we focus on such aspect. In particular, the cubature rules for bivariate weakly singular  integrals developed in \cite{Nash20} are here firstly applied to the 3D IGA-BEM context.  These rules are based on spline Quasi-Interpolation (QI) and on the usage of a spline product formula \cite{Morken91}. Specifically formulated for weakly singular integrals containing a tensor-product B-spline factor in the integrand, their peculiarity consists in being applicable on its whole support. Since  tensor product B-splines are here adopted to span the discretization space, this feature implies that the assembly phase of the system matrix can be done function-by function. This is a clear advantage of our approach, as already outlined in the 2D setting \cite{ACDS3}.
Furthermore, as a preliminary step the multiplicative singularity extraction procedure considered in \cite{Nash20} is here improved with the aim to apply the QI operator adopted by our rules to sufficiently smooth functions. Note that in the IGA-BEM setting, according to the distance between the source point and the integration domain, not all the integrals to be dealt with are singular, some being near singular and others regular. The final accuracy of our IGA-BEM implementation gained from treating near singular integrals with the same strategy adopted for singular ones. For regular integrals,  a simplified formulation of the mentioned cubature rules has been developed, relying on  the same QI operator, but possibly selecting for it a different spline degree and/or knot spacing.

The structure of the paper is as follows.
Section~\ref{sec2:indirect} recalls the indirect integral formulation of a potential problem on a 3D unbounded domain, and then it introduces the related discretization based on an isogeometric collocation BEM. Section~\ref{sec3:cub} presents cubature rules based on quasi-interpolation and on the spline product formula. Section~\ref{sec4:num} reports the results obtained with a Matlab implementation of the scheme for a screen problem and for two different  problems  on the unbounded domain external to a torus. Finally, some conclusions are given in Section~\ref{sec:conc}.

\section{Isogeometric indirect BEMs on 3D unbounded domains}
\label{sec2:indirect}
%

As well as for interior problems, even for general exterior problems different integral formulations are available in the literature, coming from different simple or double layer representation formulas.  BEMs derived from an integral equation which directly involves the unknown Cauchy datum are  called  (simple or double layer) ``direct'' methods. As an alternative, (simple or double layer) ``indirect'' approaches can also be adopted. In such case the unknown appearing in the integral equation and in the associated representation formula is a different function, usually referred to as ``density'' function. In this paper we adopt a BEM indirect formulation, since it is the only one applicable also to screen/crack problems,   \cite{costabel1986principles}. In particular we focus on such formulation for 3D Dirichlet exterior potential problems,
 \begin{equation} \label{problem}
 \left\{
 \begin{array}{l}
 \Delta u = 0\,, \,\, \x \in \RR^3 \setminus {\overline{\Omega}}\,, \cr
 u  = g\,, \,\, \x \in \Gamma\,, \cr
\left.
\begin{array}{l} \vert u(\x) \vert \rightarrow  0 \cr
 \Vert \nabla u(\x) \Vert  = O( \Vert \x \Vert^{-2})
 \end{array}
 \right\}
\, \mbox{for } \Vert \x \Vert \rightarrow \,+ \infty\,,

 \end{array}
 \right.
 \end{equation}
 where $g$ is the given Dirichlet datum and $\Gamma := \partial \Omega,$ with $\Omega$ denoting a  bounded domain in $\RR^3$.  Note that for screen problems $\Omega$ has to be replaced with  $\Gamma,$ where $\Gamma$ is a given open and limited surface. Observe also that, since we deal with unbounded domains, it is necessary to specify the required behavior of the harmonic solution at infinity  to ensure uniqueness \cite{BEMbook}.
 The boundary integral formulation of (\ref{problem}) is given by the following BIE,
 \begin{equation} \label{Vdef}
 V \psi (\x) := \int_{\Gamma} U(\x,\y) \psi(\y) \, \mathrm{d} \Gamma_{\y} \,=\, g(\x) \,, \qquad \x \in \Gamma\,,
 \end{equation}
 where the unknown density function $\psi$ represents the jump across the boundary of the normal derivative of $u$ (in the outward direction from $\RR^3 \setminus \overline{\Omega}$ to $\Omega$) and the kernel $U$ is the fundamental solution of the Laplace equation in 3D,
 \begin{equation} \label{kernel}
 U(\x\,,\,\y) := \frac{1}{4 \pi}\, \frac{1}{\Vert \x - \y \Vert_2}.
 \end{equation}
The operator $V$ defined in (\ref{Vdef}) is an elliptic isomorphism between its functional domain $H^{-\frac{1}{2}}(\Gamma)$ and its codomain $H^{\frac{1}{2}}(\Gamma),$ while the BIE also defined in (\ref{Vdef}) is a Fredholm integral equation of  the first kind with weakly singular kernel.

Once the density function $\psi$ has been (approximately) computed, the (approximated) expression of $u$ at any point in $\RR^3 \setminus \overline{\Omega}$ can be derived by relying on its simple layer representation formula,
 \begin{equation} \label{rep} u(\x) = \int_{\Gamma} U(\x \ ,\y)\, \psi(\y)\, \mathrm{d} \Gamma_{\y}\,, \qquad \x \in \RR^3 \setminus \overline{\Omega}\,.
 \end{equation}
Note that this formula automatically ensures the required decay behavior of $u$ at infinity if $\psi$ is continuous in $\Gamma$ \cite{BEMbook}, and also that its numerical approximation can be challenging only when $\x$ is very close to $\Gamma.$

To compute a numerical solution of the BIE in (\ref{Vdef}) we adopt collocation, since this is probably the easiest discretization method  suited for smooth problems (i.e., problems with a smooth $\Gamma$   and  sufficiently smooth Dirichlet datum).  However, we observe that  the related theory is not yet fully developed in particular for the 3D case;  for example in Sections 9.2 and 9.3 of \cite{Atkinson_2009} the available theoretical  results  can be applied just for the direct approach which is not usable for screen problems, as we already remarked.  In general, after having numerically solved the considered BIE with a BEM approach, the representation formula associated with it has to be considered in order  to compute values of the approximated solution inside the domain. However there are several applications, e.g. the simulation for screen or crack problems, which do not require such post-processing step, since their physically relevant quantity is the unknown Cauchy datum on the boundary itself, see for example \cite{costabel1986principles}.

We assume that the geometry is CAD-generated, that is  $\overline{\Gamma}$  has a parametric representation,
$\overline{\Gamma} = \mbox{Image}(\F)\,,$ where $\F: R  \rightarrow \RR^3$ is a one-to-one mapping\footnote{With the exception of the boundary of the parametric domain  in case of a closed $\Gamma$ or with cylinder-like shape.} defined on the parametric rectangular domain $R := [a_1\,,\,b_1] \times [a_2\,,\,b_2]$ in the following standard NURBS form,
\begin{equation} \label{NURBS} \F(\t) := \frac{ \displaystyle{\sum_{i=0}^{n_1} \sum_{j=0}^{n_2}} w_{i,j}\, \Q_{i,j}\, B^{d_1}_{i,T_1}(t_1) B^{d_2}_{j,T_2}(t_2)}{\displaystyle{\sum_{i=0}^{n_1} \sum_{j=0}^{n_2}} w_{i,j}\, B^{d_1}_{i,T_1}(t_1) B^{d_2}_{j,T_2}(t_2)}\,,\quad \t \in R\,,
\end{equation}
 where $\t :=(t_1,t_2)$. Each set of functions $B^{d_j}_{i,T_j}, i=0,\ldots,n_j, j=1,2$   is made up of  univariate B-spline functions defined in $[a_j\,,\,b_j],$ of degree $d_j$ and with extended knot vector $T_j$. We assume that the mapping $\F$ is at least differentiable.  The set $\{ \Q_{i,j} \in \RR^3, i=0,\ldots,n_1, j=0,\ldots,n_2\}$ defines the net of  control points\footnote{In order to avoid degenerate nets we assume the control net composed by distinct points.} and  each $w_{i,j}$ is a positive weight used to obtain additional shape control (obviously, if the weights are all equal to a common constant value, $\F$ becomes just a vector tensor-product polynomial spline represented in B-spline form). This surface representation  is often adopted in CAD environments to design free-form surfaces and it has the facility of involving just one patch, that is a unique continuous vector function $\F$ is used to define the whole considered surface. On this concern we note that clearly this representation can be not flexible enough for complex geometries characterized for example by linkage curves with just geometric continuity. Note also that, even if the NURBS form was originally introduced for guaranteeing exact representation of conic sections and quadric patches, it cannot  be used for example for representing closed surfaces with zero genus without singularities in the parameterization. For example in reference \cite{TauRodHug} the one-patch NURBS parameterization of the sphere is singular at the poles and this has required the development of special cubature formulas to obtain accurate results. In order to deal with surface representations more flexible and more suitable for the analysis, multi-patch NURBS representation could be considered, as well as different kinds of spline spaces, e.g., multi-degree polar splines \cite{TSHH20}.

Focusing on the basic standard single patch NURBS representation and referring for brevity only to the case of clamped knot vectors $T_j, j=1,2$, let us define ${\cal T}_j$ as a possible refinement of $T_j, j=1,2,$   obtained by inserting new breakpoints in $[a_j\,,\,b_j],$ or by increasing the multiplicity of a knot already included in $T_j.$  Denoting with $h_j$  the maximal distance
between successive knots in ${\cal T}_j, j=1,2$ and  with $S_{h_j}$  the univariate spline space of dimension $m_j+1\,, m_j \ge n_j,$ spanned by the B-splines of degree $d_j$ associated with ${\cal T}_j,$ we can define the tensor-product spline space ${\cal S}^h := S_{h_1} \otimes S_{h_2}$ of dimension $D := (m_1+1) (m_2+1),$ where   $h:= \max\{h_1,\,h_2\}.$  Note that in the following we simplify the notation by denoting as $B_j$ the bivariate B-spline $B^{d_1}_{i_1,{\cal T}_1} B^{d_2}_{i_2,{\cal T}_2}$ generating ${\cal S}^h,$ where $ j= i_2(m_1+1)+i_1.$
Thus, the density function $\psi$ is  approximated in a finite dimensional space $\hat {\cal S}^h$  derived from ${\cal S}^h$ by lifting all the B-splines to the physical space $$\hat {\cal S}^h \,:=\, \mathrm{span}\{ B_j \circ \F^{-1} , j=0,\ldots,D-1\}\,.$$
Using collocation, the coefficient vector $\bfax_h =(\alpha_0,\cdots,\alpha_{D-1})^T \in \RR^D$ associated with the IGA-BEM discrete solution $\psi_h \in \hat{\cal S}^h,\ \psi_h(\x) = \sum_{j=0}^{D-1} \alpha_j (B_j \circ \F^{-1}) (\x)\,,$ is determined by solving the following linear system with $D$ equations,
\begin{equation} \label{linsys}
A_h  \bfax_h = \g_h\,,
\end{equation}
where
\begin{equation} \label{Adef}
(A_h)_{i,j} :=   \int_{\Gamma} U(\x_i,\y) \ (B_j \circ \F^{-1}) (\y)\ d \Gamma_{\y} \,,
\end{equation}
and $(\g_h)_i := g(\x_i)\,.$

The points $\x_i , i=0,\ldots,D-1,$ are $D$ distinct collocation points  belonging to $\Gamma$ defined as the image through $\F$ of the Cartesian product of two sets of distinct abscissas, defined in $[a_1,\,b_1]$ and in $[a_2,\,b_2],$ respectively, as for example Greville abscissas (or suitable variants, see for example \cite{Wang15}, to avoid collocation on the boundary of $\Gamma$). Thus, for each collocation point $\x_i \in \Gamma, i=1,\ldots,D,$ there exists $\s_i \in R$ such that $\x_i = \F(\s_i),$ with $\s_i \ne \s_j, i \ne j.$

\noindent
Then, moving to the parametric domain, the entries of the matrix $A_h$ can be written as follows,
\begin{equation} \label{Adef2}
(A_h)_{i,j} \,=\,   \int_{R_j} U(\F(\s_i),\F(\t) )\ B_j (\t) \, J(\t) \ dt_1 dt_2\,,
\end{equation}
where $R_j := \mathrm{supp}(B_j)$ and where $J$ represents the infinitesimal surface area  element,
$$J(\cdot) := \left \Vert \frac{\partial \F}{\partial t_1}(\cdot) \times   \frac{\partial \F}{\partial t_2}(\cdot) \right \Vert_2 \,.$$
Regarding the matrix $A_h,$ we observe that it  is a full matrix whose characteristics also depend on the distribution of the collocation points in $\Gamma$. On this concern we mention that the nonsingularity of the collocation matrix descending from a second kind Fredholm BIE has been theoretically studied in \cite{Atkinson_2009}. However our $A_h$ comes from a Fredholm integral equation of the first kind and, up to our knowledge, sufficient conditions ensuring such property are not available in this case, unless for very special situations \cite{Costabel92screen3D}.  We can only add that in all our experiments we have always dealt with positive definite matrices, see also the comment about the condition number in the numerical experiments presented in Section \ref{sec4:num}.

In order to be able to compute the discrete solution efficiently and accurately, it is important to have an efficient and stable method for solving the linear system in (\ref{linsys}) but even more to consider suitable cubature rules for weakly singular integrals to be used in the assembly phase. In this paper we focus on this second issue, extending to cubature the quadrature rules previously developed for analogous univariate weakly singular integrals. As already mentioned in the introduction, the use of such rules (which are specific for integrals containing an explicit B-spline factor in the integrand) allows us to avoid the element-by-element assembly strategy, since they are directly applied on the whole support of each $B_j.$

\section{Cubature rules}
\label{sec3:cub}
%

In this section we present the  cubature rules introduced in \cite{Nash20} and here firstly applied in the 3D IGA-BEM context for the numerical approximation of all the integrals on the right hand side of (\ref{Adef2}). Such rules rely on spline quasi-interpolation and are an extension to the bivariate case of those firstly studied in the univariate setting \cite{CFSS18} and successfully tested also working in adaptive spline spaces \cite{FGKSS_2019}.

We first observe that the integrals defining the entries of $A_h$ are weakly singular only if the source point $\s_i$ belongs to ${\overline  R}_j.$ If this is not the case, they are nearly singular or regular, depending on the distance between $\s_i$ and the integration domain $R_j,$  clearly taking also into account if $\Gamma$ is a closed surface. This initial classification is important because nearly singular integrals need as much care as weakly singular ones. In our implementation we adopt for them the same strategy used for weakly singular integrals, while a simplified version of the developed rules is used to approximate regular integrals, possibly using also higher QI spline degree in such a case.

Referring first of all to singular and nearly singular integrals, we observe that, as well as in the 2D case \cite{ACDS3, FGKSS_2019}, a preliminary singularity extraction procedure is necessary in order to increase the performance of our cubature rule. In \cite{Nash20} an analysis has been developed to show that two singularity extraction procedures,  based on a Taylor expansion of $\F(\t)$ about $\s$, are possible, one using an additive splitting of $U$ and the other of multiplicative type. Considering that the first splitting requires the numerical computation of two integrals instead of one, here we have chosen to rely on the multiplicative technique. In this case the idea is to rewrite the entries of the matrix $A_h$ given in (\ref{Adef2}) in the following equivalent form,
\begin{equation} \label{Adefnew}
(A_h)_{i,j} :=   \int_{R_j} K_\F(\s_i,\t)\  \rho_{\s_i,\F}(\t) B_j (\t) \, J(\t) \ \mathrm{d}t_1 \mathrm{d}t_2\,,
\end{equation}
where
$$\rho_{\s,\F}(\t) := \frac{U(\F(\s),\F(\t))}{K_\F(\s,\t)}\,,$$
and $K_\F$ is a suitably locally defined weakly singular kernel such that the function $\rho_{\s,\F}$ is at least continuous.
In particular in \cite{Nash20} we considered the following definition of such reference kernel,
\begin{equation} \label{KFdef}
K_\F := K_\F(\s,\t) :=  \frac{1}{4\pi} \ \frac{1}{\sqrt{P_{\s,\F}(\t)}}\,,
\end{equation}
with $P_{\s,\F}$ denoting the bivariate homogeneous quadratic polynomial vanishing with its first derivatives at $\t = \s$ and obtained by replacing $\F(\t)$ with its first degree Taylor expansion about $\s$ in $\Vert \F(\t) - \F(\s)\Vert_2^2$. Note anyway that for a general surface this definition just ensures continuity to $\rho_{\s,\F}$ at $\t=\s$. However, if $\F$ is sufficiently smooth at the singular point $\s$, it is possible to use higher order expansions for $\F$, in order to ensure higher regularity to $\rho_{\s,\F}$. In particular here we consider an expansion with up to $3$ terms, which ensures up to $C^2$ smoothness to $\rho_{\s,\F}$. Details about this improved singularity extraction procedure will be presented in a forthcoming paper\footnote{T. Kandu\v{c}. Isoparametric singularity extraction technique for 3D potential problems in BEM, arXiv:2203.11538.}, since they form a technical part necessary for an exhaustive explanation but too long for being here reported. We just mention that in this case
$K_\F$ is defined with the following expression which reduces to (\ref{KFdef}) for $n=1,$
\begin{equation} \label{newKF}
K_\F(\s,\t) :=  \frac{1}{4\pi} \ \sum_{\ell=1} ^n \sqrt{P_{\s,\F}(\t)}^{-2 \ell + 1} P_{3\ell-3}^{(\s,\F)}(\t),
\end{equation}
where $P_{3\ell-3}^{(\s,\F)}$ are appropriate homogeneous polynomials of total degree $3\ell-3\,, \ell=1,\ldots,n.$
Concerning this improved singularity extraction procedure, finally we add that, when $n >1 $ and the computational mesh is not fine enough,  a correction term $ \eta \|\s - \t \|_2^2 $ for some suitably small $\eta>0$ has to be added on the right hand side of (\ref{newKF}) to avoid the introduction of new singularities in $\rho_{\s,\F}$.

After this preliminary rewriting of the entries of the matrix $A_h,$   each of them is approximated as follows,
\begin{equation} \label{entryapprox}
(A_h)_{i,j} \approx   \int_{R_j} K_\F(\s_i,\t)\  \left( Q_{\p,\Delta_j} \left(J  \rho_{\s_i,\F} \right)\right) (\t)\ B_j (\t) \ \mathrm{d}t_1 \mathrm{d}t_2\,,
\end{equation}
where $\Delta_j$ is a uniform partition of $R_j$ and $Q_{\p,\Delta_j}$ denotes a tensor-product  QI operator which approximates a bivariate function  with a  tensor-product spline of bi-degree $\p = (p_1,p_2)$  with respect to the  partition $\Delta_j.$

 In order to compute the integral on the right hand side of (\ref{entryapprox}), the spline $Q_{\p,\Delta_j } \left(J \  \rho_{\s_i,\F} \right)$ is  multiplied by the B-spline $B_j,$ using for this aim the extension to the tensor-product setting of the M{\o}rken formula \cite{Morken91} for spline product. As a result, the whole product $\left(J \  B_j\ \rho_{\s_i,\F}\right) $ is approximated with a local tensor-product spline function defined in $R_j,$ having bi-degree $\p + \d,$ and breakpoints in $R_j.$
We observe that, in order to implement the so derived cubature rules, the evaluation of the following {\it modified moments} is firstly necessary,
\begin{equation}  \label{MM}
\mu_{K_\F,B}(\s) := \int_{ \mathrm{supp}(B)} K_\F(\s,\t) \ B(\t) \ \mathrm{d}t_1 \mathrm{d}t_2 \,,
\end{equation}
with $B$ denoting any tensor-product B-spline. This computation is done by relying on the recursion formula for B-splines and on the preliminary analytic computation of element integrals whose integrand is the product of a monomial times $K_\F,$ procedure possible for a kernel $K_\F$ defined as in (\ref{KFdef}) or more generally as in (\ref{newKF}).

When the integral defining the $(i,j)$--th entry of $A_h$ is classified as a regular one,  we do not need the preliminary computation of the modified moments and the strategy is simplified by directly setting
\begin{equation} \label{entryappr_reg}
(A_h)_{i,j} \approx   \int_{R_j}  \left( Q_{\q,\Lambda_j} \left(U(\F(\s_i),\cdot) J  \right)\right) (\t)\ B_j (\t) \ \mathrm{d}t_1 \mathrm{d}t_2\,,
\end{equation}
where $\q =(q_1\,,\,q_2)$ denotes a  bi-degree for the local splines possibly different from $\p$ and $\Lambda_j$ is a uniform partition of $R_j$, also possibly different from $\Delta_j.$

Concerning the QI spline operator used by the proposed cubature rules, for our experiments we have always considered the uniform tensor-product formulation of the derivative free variant of the Hermite QI operator introduced in \cite{MSbit09}, see for example \cite{BGMS16} for an introduction to its formulation in the bivariate setting and for an application in the context of adaptive bivariate spline approximation.
This univariate spline QI operator  has approximation order $p+1$ and in  the context of numerical integration it has special interest because its application to both regular and weakly singular quadrature is superconvergent when even spline degrees are considered on uniform knot distributions \cite{MSJcam12,CFSS18}. Thus for smooth functions $Q_{\p,\Delta_j}$ has approximation order $\min\{p_1,p_2\}+1$  and it inherits the superconvergence feature when analogously applied for cubature. Observe however that, in our approach, choosing values of $p_k, k=1,2,$ greater than $2$ is not useful, because of the limited regularity of the function $\rho_{\s_i,\F}.$  It can be profitable instead to select $q_k >2, k=1,2$, if $\F$ is smooth in $R_j.$ Considering the assumed analytic expression of $\F,$ this is surely true when the interior of $R_j$ does not contain any straight line $t_j=\tau_j,$ with $\tau_j$ {\it geometric knot} in the $j$--th parametric direction, that is $\tau_j \in T_j\,, j=1,2,$ where $T_j, j=1,2$ are the extended knot vectors used to define $\F.$

\section{Numerical simulations}
\label{sec4:num}

In this section we test the uniform formulation of the proposed IGA boundary element model for the numerical solution of two Laplace problems on two different 3D unbounded domains. Furthermore a nonuniform implementation of the scheme is also tested for solving a more difficult Laplace problem exterior to the first considered domain. In all the examples the improved singularity extraction procedure is used within the developed singular cubature rules, always setting the parameter $\eta$ (see Section \ref{sec3:cub}) equal to $1.$ Concerning the discretization, in all the presented experiments approximate densities are constructed on a given initial mesh and also on successive finer ones obtained from a dyadic $h$-refinement procedure which inserts a simple knot at the midpoint between every pair of consecutive breakpoints in each parametric direction. \\
For the first two considered examples the density $\psi$ is known, while the Dirichlet datum $g$ is not available in closed form. Thus, the entries of the right-hand side in \eqref{linsys} are computed numerically,
\begin{equation}
(\g_h)_{i} =   \int_{R} U(\F(\s_i),\F(\t))\ (\psi \circ \F)(\t)\, J(\t) \ \mathrm{d}t_1 \mathrm{d}t_2\,,
\end{equation}
 by applying the same approach adopted to compute the system matrix. For these examples the accuracy of the scheme is measured by the $L^2(\Gamma)$ norm of the absolute density error $(\psi_h - \psi),$ in particular considering its dependency on the mesh size $h.$  The third exterior Laplace problem is based on Test $5.3$ of \cite{Dolz18} where  the more general Helmholtz equation is solved. Since in this case the analytic expression of the exact potential $u$ is known, the Dirichlet datum is given and the accuracy in the potential approximation $u_h \approx u$ can be measured, provided that the representation formula given in (\ref{rep}) is previously applied to $\psi_h$ to define $u_h$ in  the considered unbounded domain $\RR^3 \setminus \overline{\Omega}.$\\

.\\\subsection{Problem on an unbounded domain exterior to a torus}\label{ex1}

\begin{figure}[t!]
\centering
\subfigure[Initial mesh in the parametric domain $R$]{
{\includegraphics[trim = 1cm 1.25cm 1cm 1.25cm, clip = true, height=3.5cm]{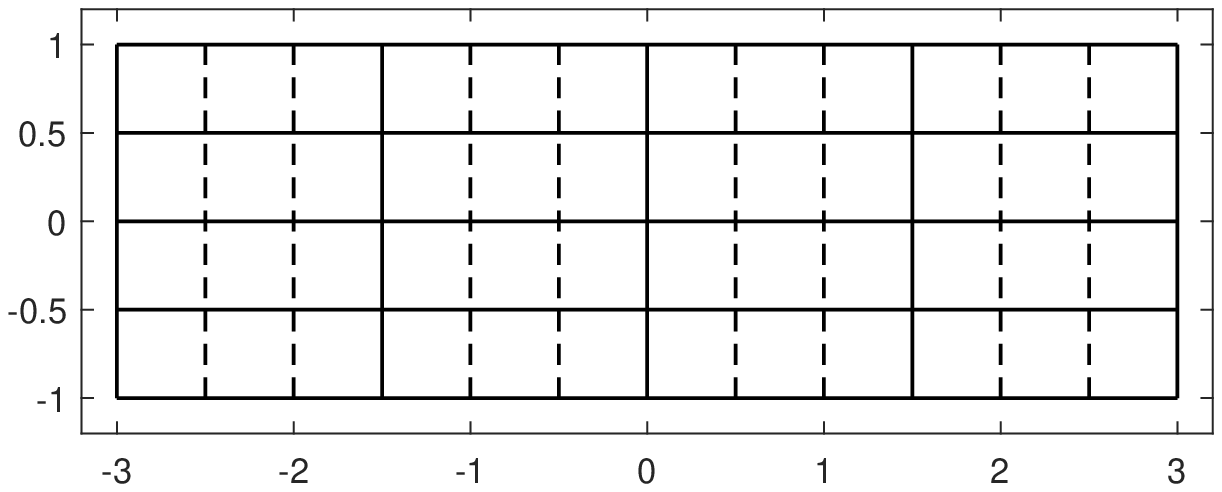}}
}
\subfigure[Control net and weights for one patch]{
{\begin{overpic}[trim =1cm .5cm .5cm 2cm, clip = true, height=3.5cm]{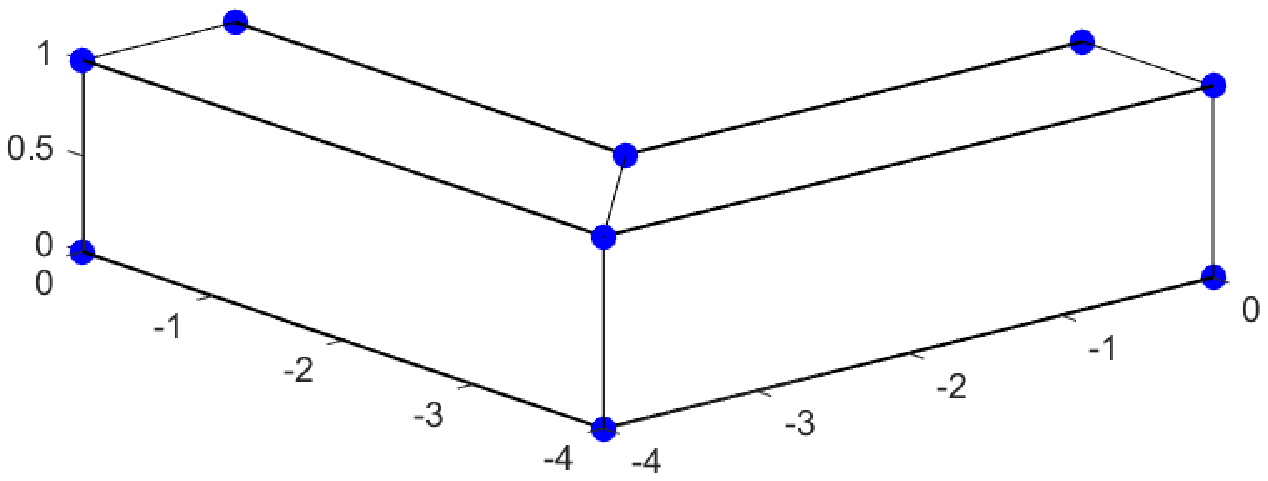}
\put(12,41){\scriptsize 1}
\put(6,77){\scriptsize $\sqrt{2}/2$}
\put(30,80){\scriptsize 1}

\put(65,18){\scriptsize  $\sqrt{2}/2$}
\put(75,48){\scriptsize $1/2$}
\put(72,67){\scriptsize  $\sqrt{2}/2$}

\put(140,37){\scriptsize  $1$}
\put(127,61){\scriptsize $\sqrt{2}/2$}
\put(123,76){\scriptsize  $1$}
\end{overpic}}
}
\subfigure[Surface $\Gamma = \partial \Omega = \F(R)$ and initial mesh]{
{\includegraphics[trim = 1cm 1.25cm .5cm 1.25cm, clip = true, height=3.5cm]{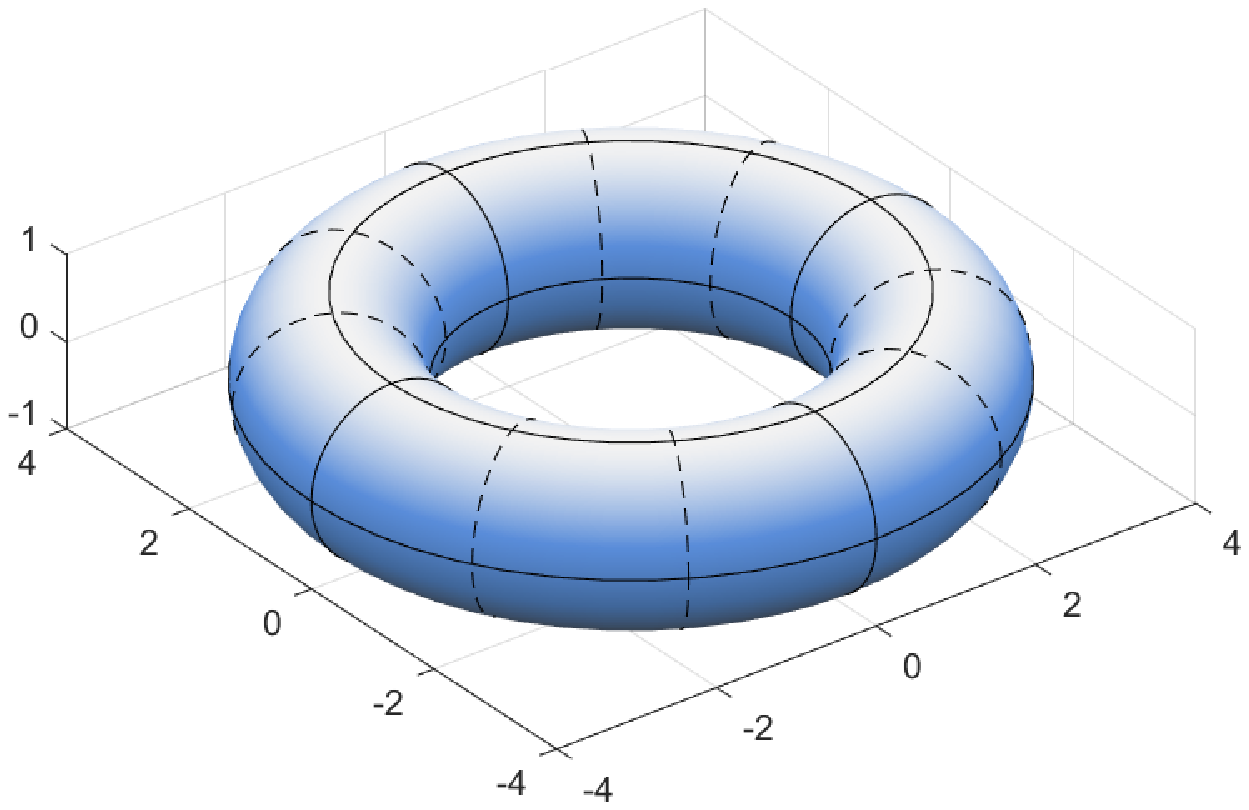}}
}
\subfigure[Graph of composition with $\F$ of the exact BIE solution $\psi$]{
{\includegraphics[trim = 0cm .0cm 0cm 0.0cm, clip = true, height=3.5cm]{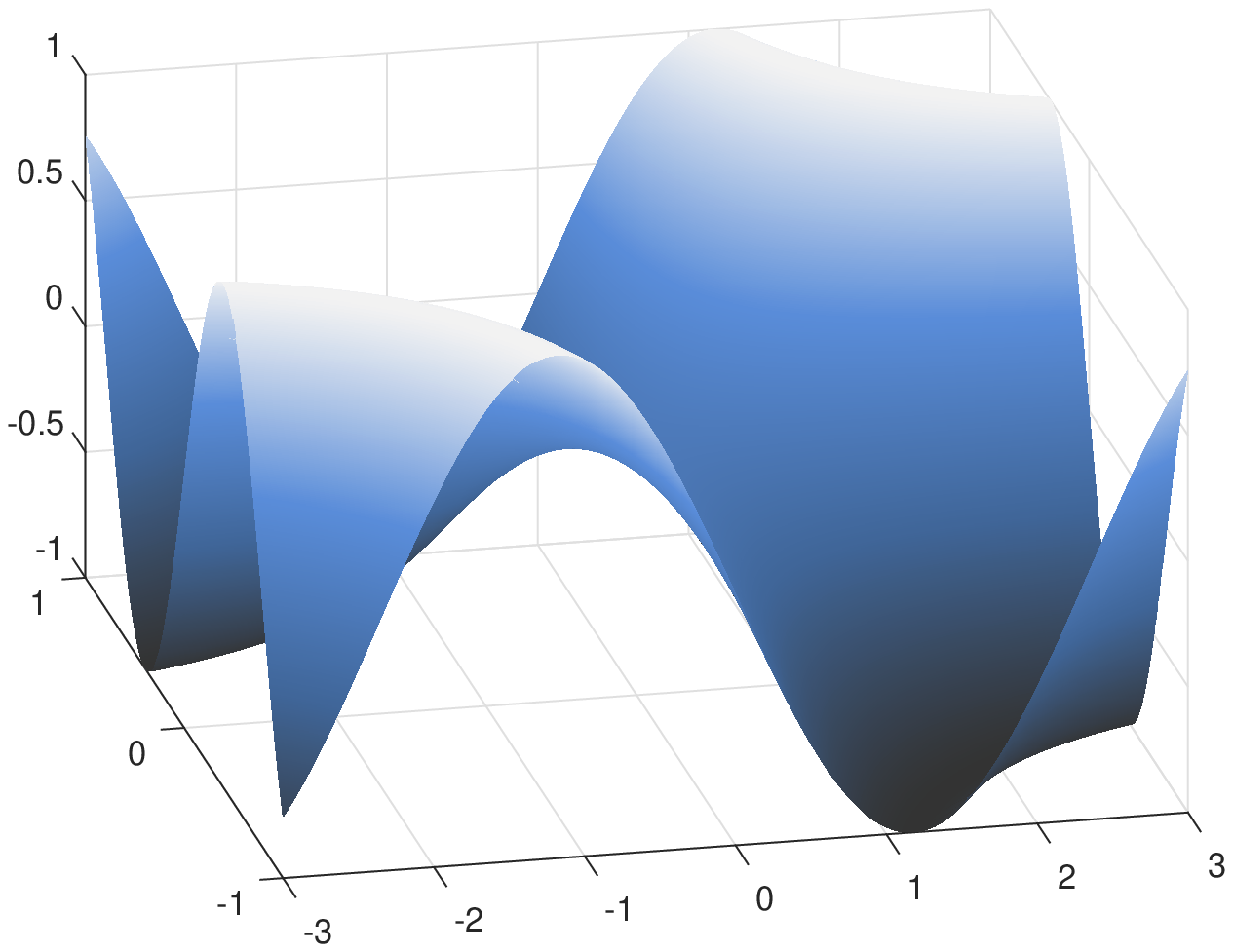}}
}
\caption{Example \ref{ex1}: problem exterior to a torus.}
\label{fig:torus}
\end{figure}

We consider a problem exterior to a toroidal surface $\Gamma$ with a major radius $\rho_M = 3$ and a minor radius $\rho_m = 1$. The geometry is represented in NURBS form by defining as in (\ref{NURBS}) the vector function $\F$  on the parametric domain $R = [-3\,,\,3] \times [-1\,,\,1]$ with  $\d = (2,2),$ fixing the following knot vectors,
\begin{align*}
T_1 = (   -3,\;   -3,\;   -3,\;   -1.5,\;   -1.5,\;   -1.5,\;    0,\;         0,\;         0,\;     1.5,\;    1.5,\;    1.5,\;    3,\;    3,\;    3),\\
T_2 = ( -1,\;   -1,\;   -1,\;   -0.5,\;   -0.5,\;   -0.5,\;         0,\;         0,\;         0,\;    0.5,\;    0.5,\;    0.5,\;    1,\;    1,\;  1),
\end{align*}
see solid lines in Fig.~\ref{fig:torus}(a).
In this form the torus is composed by the $16$ rational biquadratic sub-patches shown in Fig.~\ref{fig:torus}(c) (separated by solid lines), which are locally represented in rational Bernstein form, since all the breakpoints have maximal multiplicity (equal to $3$) in both the knot vectors $T_j, j=1,2.$
The associated control net has the size $12 \times 12;$  for clarity in Fig.~\ref{fig:torus}(b) we show a sub-patch and the related $3 \times 3$ sub-net of control points with the corresponding weights. Dashed lines in Fig.~\ref{fig:torus}(a) and Fig.~\ref{fig:torus}(c) represent additional inserted knots for the initial computational mesh. Note that, despite the maximal multiplicity of the breakpoints in $T_j, j=1,2,$ the vector function $\F$ is $C^1$ continuous because of the symmetry of the considered control net and of the related weights. We remark that, even if not strictly necessary, maximal multiplicity is used for the geometric knots to ensure to $\F$ maximal regularity when restricted to the support of a trial function $B_j.$ As explained at the end of the previous section, this is useful for the accuracy of our cubature rules when  regular integrals are dealt with, since the support of a trial function can never cross a line $t_j = \tau_j \in T_j, j=1,2\,,$ i.e. solid lines in Fig.~\ref{fig:torus}(a). \\
The exact function $\psi \circ \F$ is set to
\begin{align*}
(\psi \circ \F)(\t) = \sin(4/3\; t_1\, t_2 )
\end{align*}
and its shape in the parametric domain is depicted in Fig.~\ref{fig:torus}(d). Each rational sub-patch of the torus
contains the same number of collocation points obtained from the Cartesian product of the improved Greville abscissas \cite{Wang15}, defined on the corresponding sub-rectangles of $R$, and mapped to $\Gamma$ by $\F.$\\
Two different experiments have been performed for this example. The first one is a standard test where initially ${\cal T}_2 = T_2$ and ${\cal T}_1$ is obtained by merging $T_1$ with $\{-2.5\,, -2\,, -1\,, -0.5\,, 0.5\,, 1\,, 2\,, 2.5\}.$ Since in all the simulations the geometric knots keep  maximal multiplicity, at any refinement level also discontinuous functions belong to the space ${\cal S}_h.$ In the second experiment initially we define  ${\cal T}_2$ still from $T_2$ but by decreasing the multiplicity of each inner knot by one. ${\cal T}_1$ is initialised as in the other case but again by reducing by one the multiplicity of all the geometric knots. Thus  $3 \times 3$ degrees of freedom are eliminated  with respect to the previous experiment and as a consequence the simulation space is always included in $C^0(R).$ Note that anyway we keep unchanged the set of collocation points so we end up with an overdetermined linear system which is solved in the least squares sense.  

In both experiments the results shown are obtained by deriving the reference kernel $K_\F$ by using the first two terms of the Taylor expansion of $\F.$ For the weakly or near singular cubature the QI bi-degree $\p$ is set to $\p = (2,2),$ using in each parametric direction 13 uniform nodes on the support of any B-spline trial function.  Due to higher global smoothness of integrands, for regular integrals we can profitably set $\q = (4,4)$ and just use 7 uniform nodes on the support of any B-spline in each parametric direction.

Convergence plots of the behaviour of approximation error $\|\psi - \psi_h\|_{L^2(\Gamma)}$ using uniform $h$-refinement strategy together with the  theoretical convergence order 3 is shown in Fig.~\ref{fig:convergence}(a). Looking at the figure, we see that both the plots have the expected behaviour and also that the $C^0$ one is very near to the other.
The system matrix $A$ is well conditioned for all tests (in the last step its condition number for the first experiment is about $5.1 \cdot 10^3$).

\subsection{A screen problem for a curved obstacle}\label{ex2}
\begin{figure}[t!]
\centering
\subfigure[Surface $\Gamma = \partial \Omega = \F(R)$]{
{\includegraphics[trim = 0cm .0cm 0cm 0.0cm, clip = true, height=3.5cm]{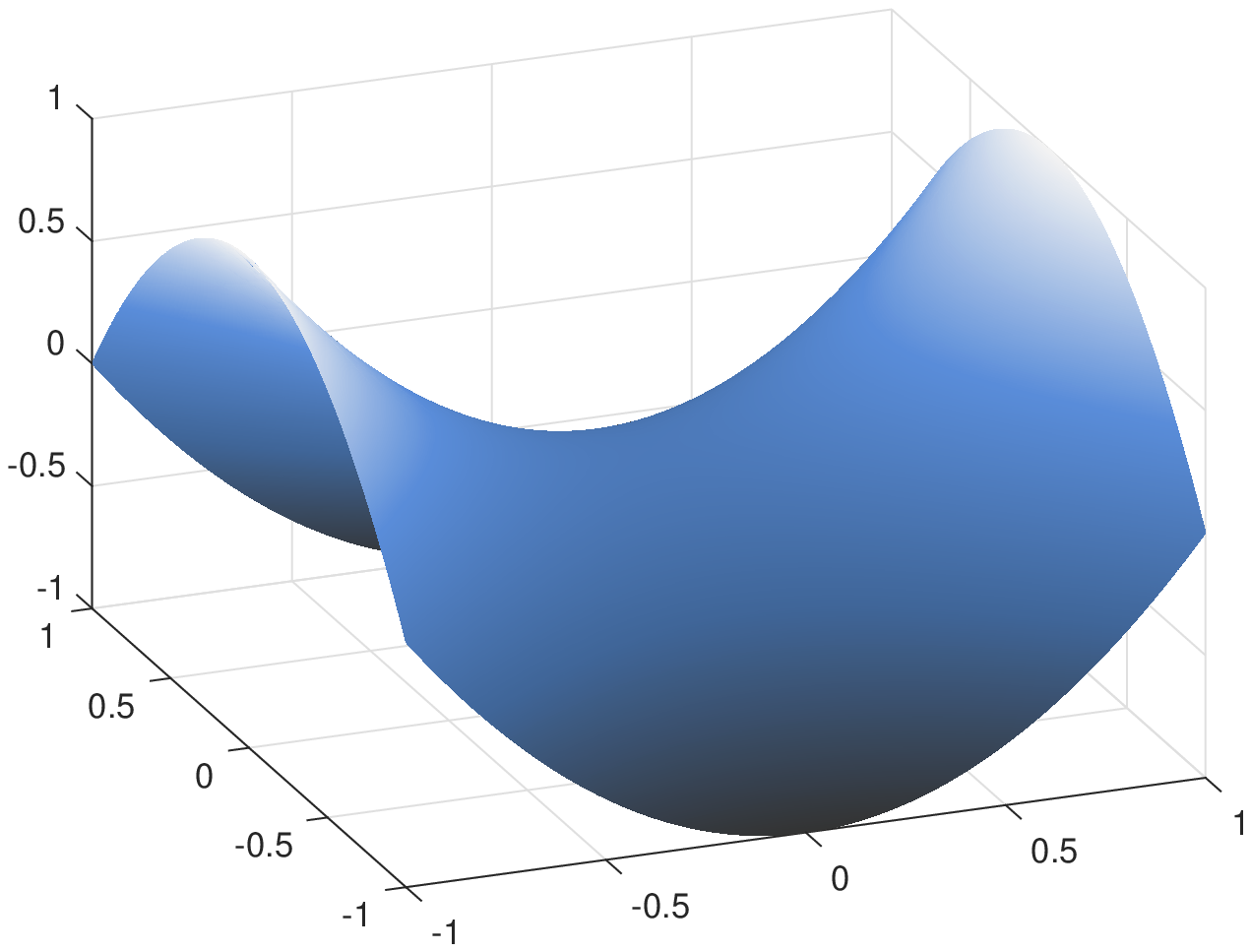}}
}
\subfigure[Graph of composition with $\F$ of the exact BIE solution $\psi$]{
{\includegraphics[trim = 0cm .0cm 0cm 0.0cm, clip = true, height=3.5cm]{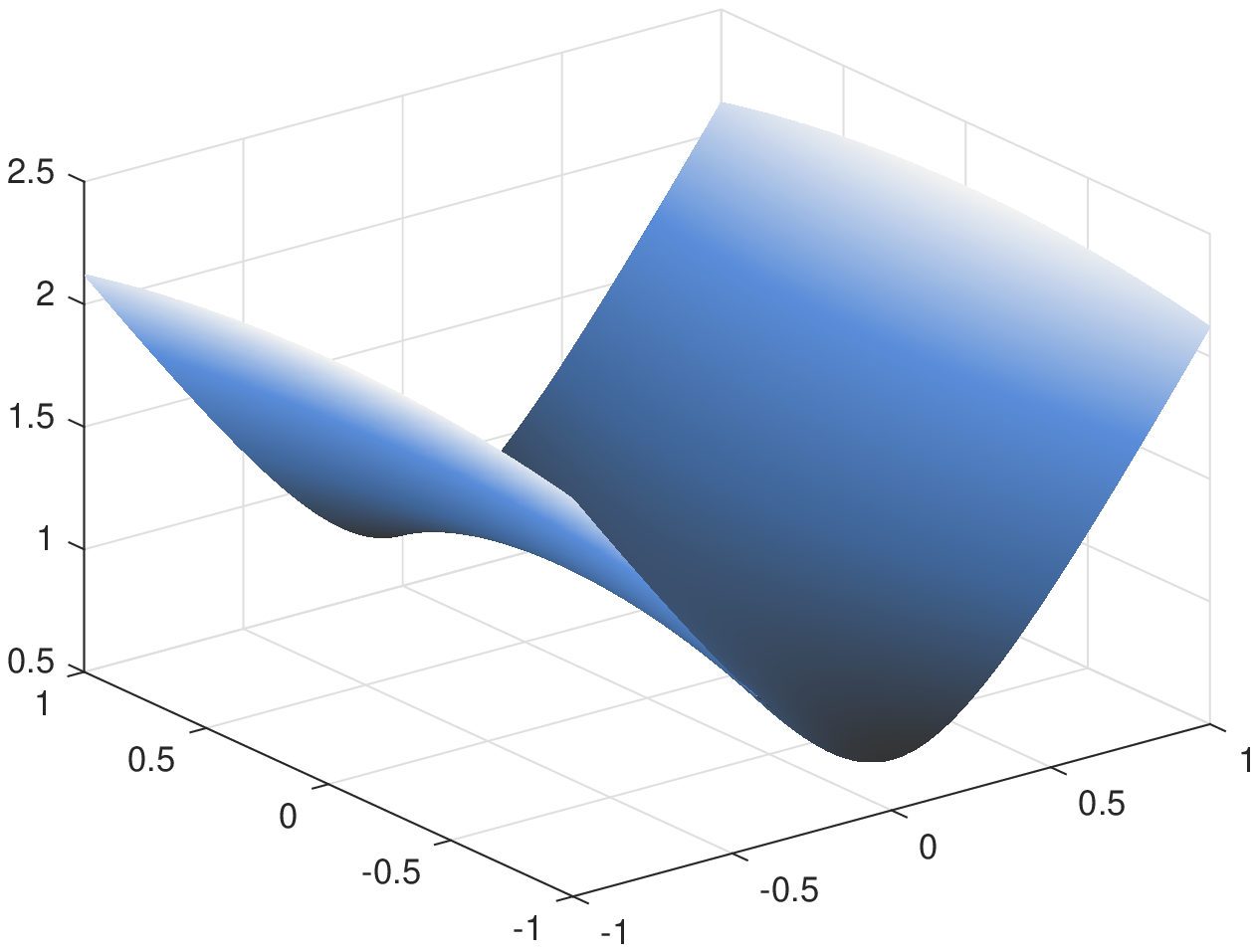}}
}
\caption{Example \ref{ex2}: screen problem.}
\label{fig:saddle}
\end{figure}

As a second example we consider a screen problem \eqref{problem} with $\Gamma$ taken as the limited portion of a saddle open surface parametrized by $\F : R \rightarrow \overline{\Gamma}$, with $R = [-1\,,\,1]^2$ and $\F(t_1,t_2) = (t_1\,,\, t_2\,,\, t_1^2 - t_2^2)$. The Dirichlet datum $g$ is numerically computed such that the composition with $\F$ of the exact density function, solution of Eq. \eqref{Vdef}, is given by
$$
(\psi \circ \F)(\t) = \sqrt{1+4t_1^2-0.5t_2^2}.
$$
 In Fig.~\ref{fig:saddle} we show the geometry and   $\psi \circ \F.$  Note that the considered $\F$ can be written as in  (\ref{NURBS}) with $\d = (2,2)$ and $T_j = (-1,\, -1,\, -1,\, 1,\, 1,\, 1)\,,\, j=1,2\,,$ since polynomials (restricted to $R$) are a subset of NURBS. Thus the problem is discretised by using quadratic B-splines and setting $\p = (2,2)$ for the QI approximation of singular and nearly singular integrals, with $7$ quadrature uniform nodes on the support of any B-spline trial function in each parametric direction.
The reference kernel $K_\F$ is obtained in this case using the first three terms in the Taylor expansion of $\F,$ while the improved  Greville abscissas \cite{Wang15} are employed as collocation points also in this case.

\begin{figure}[t!]
\centering
\subfigure[Problem exterior to a torus $\big(\d = (2,2)\big)$]{
{\includegraphics[trim = 0cm .0cm 0cm 0.0cm, clip = true, height=3.75cm]{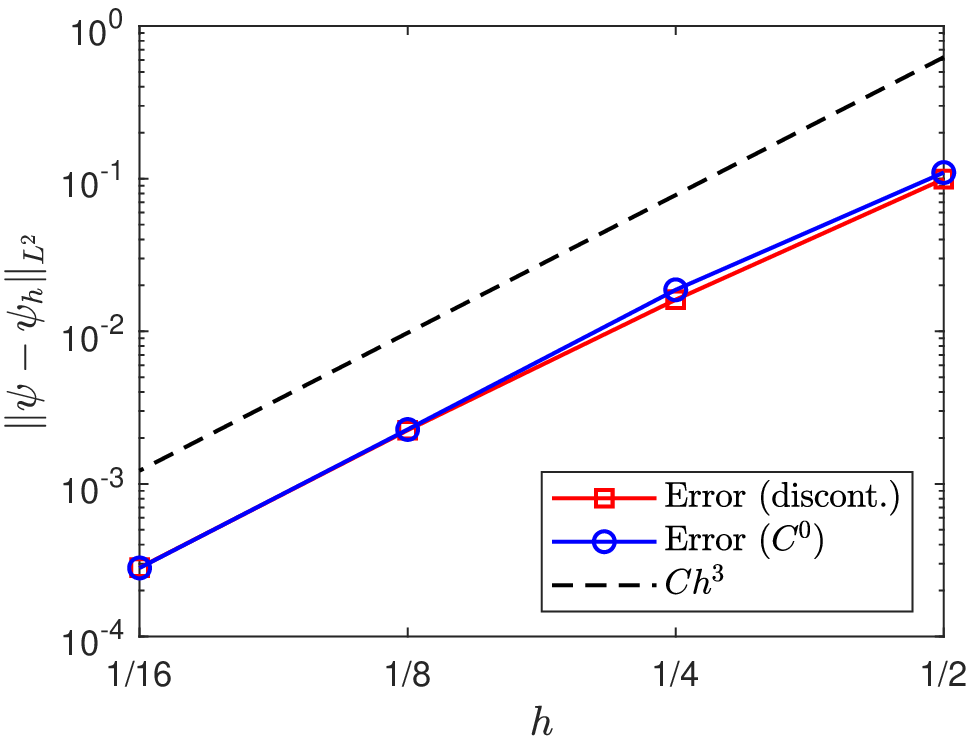}}
}
\subfigure[Screen problem $\big(\d = (2,2)\big)$]{
{\includegraphics[trim = 0cm .0cm 0cm 0.0cm, clip = true, height=3.75cm]{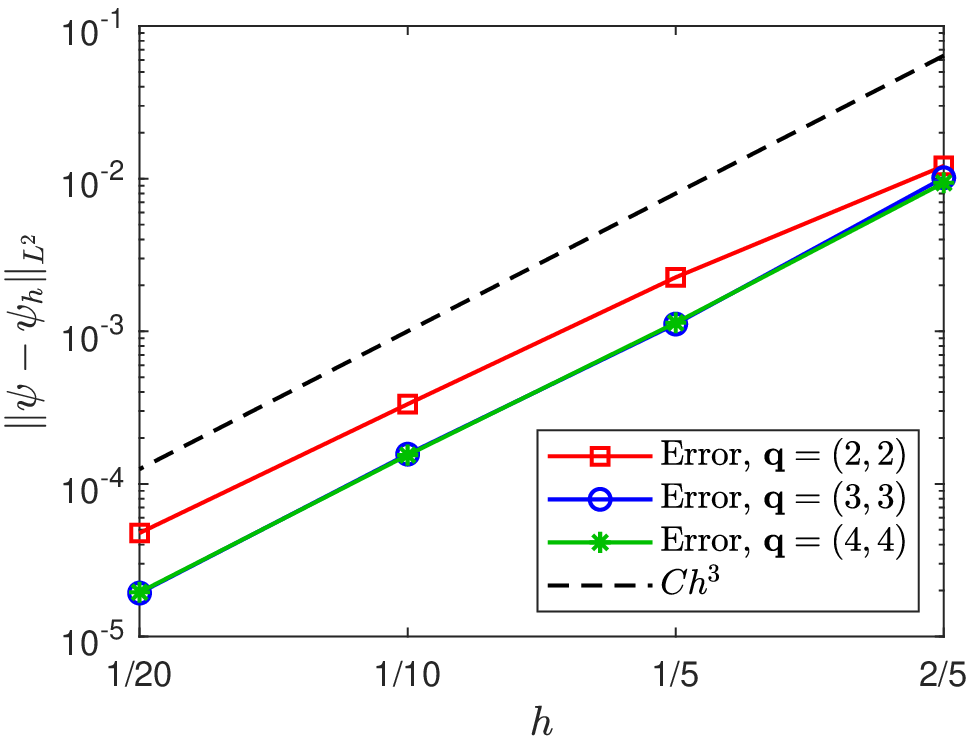}}
}
\caption{Convergence plots in $L^2(\Gamma)$ norm of the density error $\psi - \psi_h$ for Examples \ref{ex1} and \ref{ex2}.}
\label{fig:convergence}
\end{figure}

 In Fig. \ref{fig:convergence}(b) we show the convergence plots for $\q$ varying from $(2,2)$ to $(4,4)$, using from $5$ to $7$ uniform quadrature nodes per B-spline support in each parametric direction for regular integrals. Better accuracy of the numerical solution corresponds to higher $\q$ but the order $3$ is observed in all the cases. Finally, the condition number for the system matrix $A$ has order $10^2$ at the last level.

\subsection{A nearly singular potential problem exterior to a torus}\label{ex3}

\begin{figure}[h!]
\centering
\subfigure[Composition with $\F$ of $g$]{
{\includegraphics[trim = 0cm .0cm 0cm 0.0cm, clip = true, height=3.75cm]{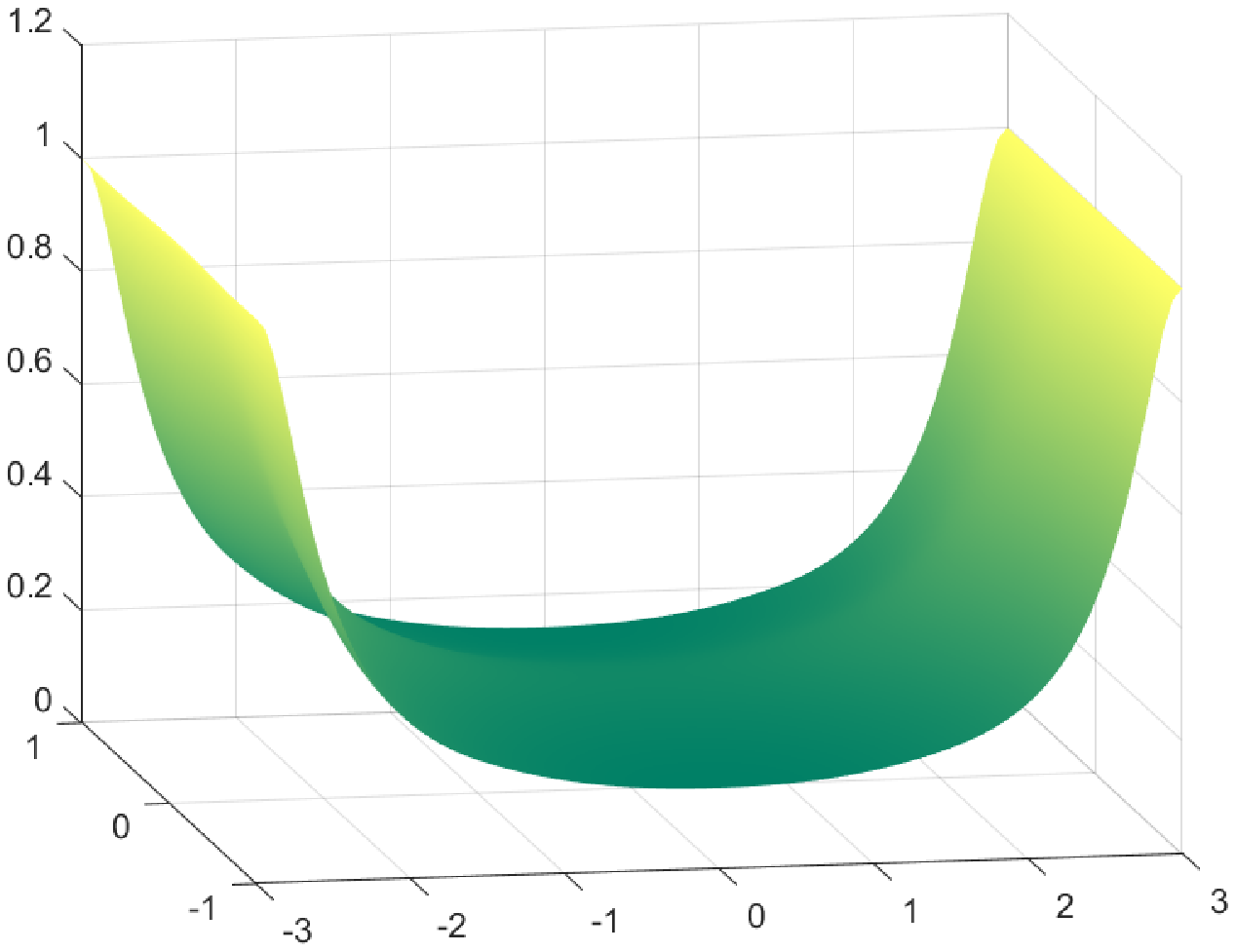}}
}
\subfigure[Composition with $\F$ of  $\psi_h$ (mesh from three refinements
of the initial one)]{
{\includegraphics[trim = 0cm .0cm 0cm 0.0cm, clip = true, height=3.75cm]{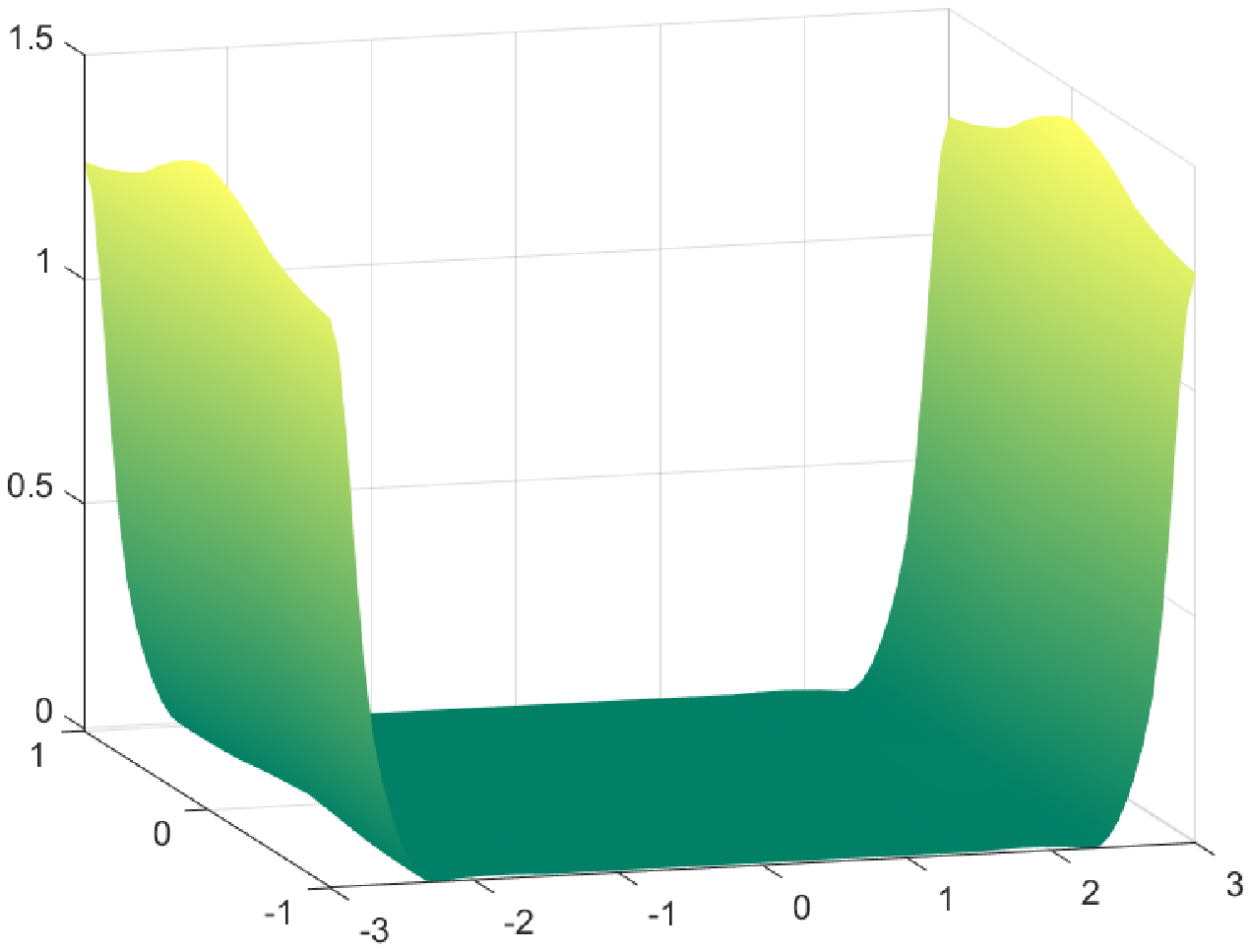}}
}
\subfigure[Initial mesh in the parameter domain]{
{\includegraphics[trim = 0cm .0cm 0cm 0.0cm, clip = true, height=3.75cm]{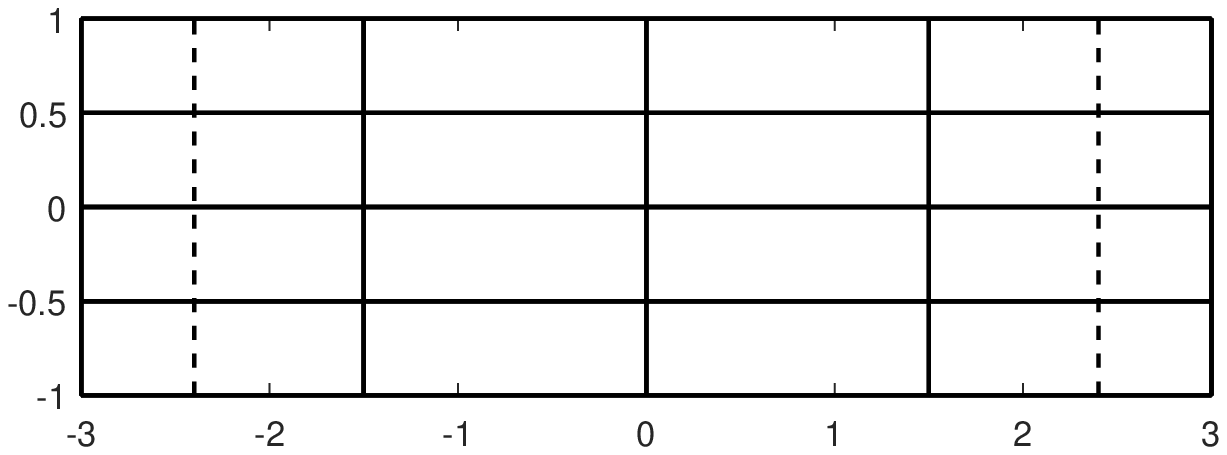}}
}
\subfigure[Distribution on $S^2$ of the potential error $\vert u - u_h \vert $ (mesh from three refinements
of the initial one)]{
{\includegraphics[trim = 0cm .0cm 0cm 0.0cm, clip = true,height=3.75cm]{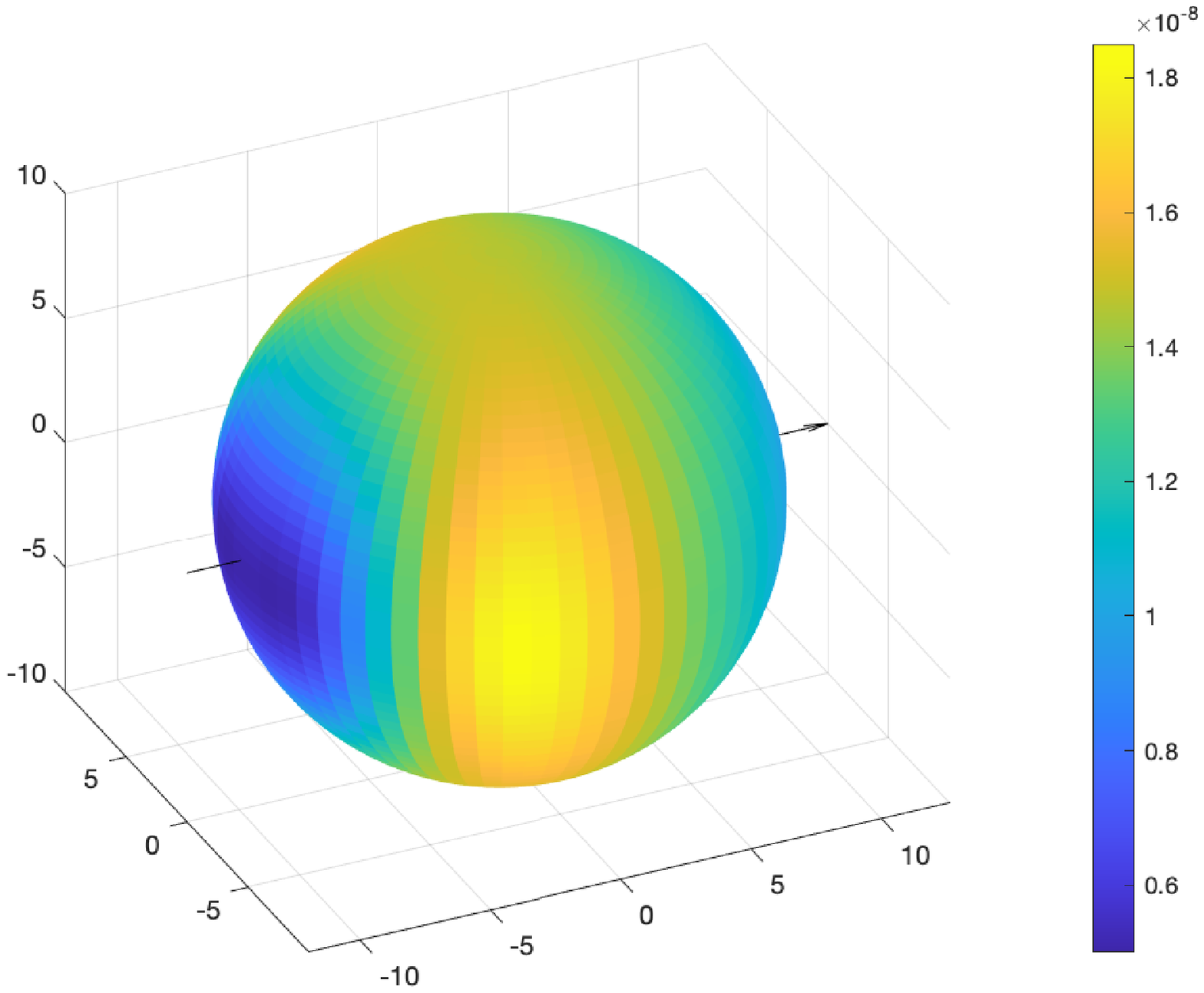}}
}
\caption{Example \ref{ex3}: another problem exterior to the torus considered in Example \ref{ex1}.}
\label{fig:layers}
\end{figure}

We consider the same toroidal surface presented in Example \ref{ex1},  solving the problem \eqref{problem} with the Dirichlet datum $g$ chosen as,
\begin{align}\label{diri-dat}
g(\x) :=  \frac{1}{\Vert \x -\bnu \Vert_2}\quad\mbox{for }\,\, \x \in \Gamma\,,
\end{align}
which is also the known analytic expression of the exact potential $u$  when $\x \in \RR^3 \setminus {\overline{\Omega}}.$ In particular in the reported experiments we fix $\bnu = (-3,0,0) \in \Omega\,,$ see Fig.~\ref{fig:layers}(a) for the graph of $g \circ \F$ with this setting. This type of Dirichlet datum gives rise to sharp gradients of the corresponding function $\psi \circ \F$ near two opposite boundary edges of the parametric domain, see  Fig.~\ref{fig:layers}(b) where we can only show the graph of $\psi_h \circ \F,$ since the exact expression of $\psi$ is unknown.   

A suitable mesh, which is nonuniform in the horizontal parametric direction, has been manually constructed in order to get better accuracy. Clearly  it would be advisable to rely on a fully automatic adaptive mesh generation strategy  but this is out of the scope of the presented paper. The initial mesh used for the reported experiments is shown in Fig.~\ref{fig:layers}(c); as in Example \ref{ex1}, solid lines correspond to the knots of the vectors $T_1$ and $T_2$ which are defined as before, while  dashed lines correspond to additional inserted knots. As done in the second subtest described in Subsection \ref{ex1}, the multiplicity of each inner geometry knot is decreased by one in order to deal with continuous spline spaces.  The collocation points are
obtained from the improved Greville abscissas and their number matches the global number of degrees of freedom characterising  the discretization space.  

Using the representation formula \eqref{rep}, for this test we also evaluate the approximated potential at some point $\x$ of the unbounded domain,
\begin{align}\label{newRep}
u_h(\x) :=  V \psi_h (\x)\quad \mbox{for }\x \in \RR^3 \setminus {\overline{\Omega}}\,.
\end{align}
In particular, in order to evaluate a posteriori the accuracy of our scheme, we compute the potential error $u -u_h$   at $2^{m+2}$ points uniformly sampled on a spherical surface $S^2$ with center located at the origin and with radius equal to $10$ and so exterior to the considered toroidal surface. Note that here $m$ denotes the dyadical refinement step, where $m = 1$ for the initial mesh and in particular  $m = 4$ was adopted to get the error distribution in $S^2$ shown in Fig.~\ref{fig:layers}(d). Looking at the figure, it can be observed that the potential error has a nice constant order of magnitude on the sphere equal to $10^{-8}$ and that its distribution is symmetric with respect to the $x$ axis where the point $\bnu$ is located (it is not shown in the figure since it is inside the sphere).

The   accuracy achieved at different refinement levels $m$ is reported in Table~\ref{tab:accuracy} together with the corresponding  global (Ndof) and directional (Ndof$_1$, Ndof$_2$) number of degrees of freedom. It is measured in both the maximum and the $L^2(S^2)$ norms, both approximated using the available discrete information on the sphere. The test has been performed by using quadratic discretization splines and selecting $\p=\q=(2,2)$.  Concerning the cubature rules,  $7$ and $13$ uniform nodes on the support of every  basis function are used in each parametric direction, respectively for the evaluation of singular/nearly singular and regular integrals. It is worth noticing that, since the evaluation of $u_h$ is done at points located far enough from the surface $\Gamma$, the integral in \eqref{newRep} is non singular. Hence, the same cubature rule constructed for the evaluation of the regular integrals has been employed for the representation formula. 
Finally for completeness we mention that also for this example the coefficient matrix in (\ref{linsys}) is well conditioned in all the developed numerical experiments, having a condition number about equal to $3.3 \cdot 10^3$ at the last considered level.
\begin{table}[h]
\begin{center}
\begin{tabular}{|c||c|c|c|c|c|}
\hline
$m$ & Ndof$_1$ & Ndof$_2$ & Ndof & $\|u-u_h\|_{L^2(S^2)}$& $\|u-u_h\|_{\infty}$  \\
\hline
$1$ &$11$&	$9$&	$99$ &$1.24\cdot10^{-2}$ & $5.88\cdot 10^{-4} $\\
\hline
$2$ & $17$        &	$13$& $221$	&$3.58\cdot 10^{-3}$& $2.21\cdot 10^{-4} $\\
\hline
$3$ & $29$& 	$21$& $609$ &$1.51\cdot 10^{-4}$& $1.00\cdot 10^{-5} $\\
\hline
$4$ & $53$ &	$37$& $1961 $ &$4.68\cdot 10^{-7}$& $1.82\cdot 10^{-8} $ \\
\hline
\end{tabular}
\end{center}
\caption{Example \ref{ex3}: degrees of freedom and norms of the potential error 
$u - u_h$ obtained at different refinement level $m$ setting $\d = \p = \q = (2,2).$} 
\label{tab:accuracy}
\end{table}


\section{Conclusion} \label{sec:conc}
In this paper a collocation IGA-BEM scheme relying on the indirect integral formulation of Laplace problems is proposed and tested for 3D exterior problems, obtaining the expected convergence order in the density approximation.
Cubature rules based on spline quasi-interpolation and on a spline product formula are applied for the numerical approximation of both weakly and nearly singular bivariate integrals appearing in the 3D IGA-BEM scheme. They are combined with a robust multiplicative singularity extraction procedure to increase the accuracy of the rules. A variant of such rules is proposed for the regular integrals to be also dealt with in the assembly phase. This can be efficiently implemented in a function-by-function form thanks to the peculiarity of the considered cubature.


\begin{acknowledgement}
A. Falini, C. Giannelli, M. L. Sampoli, and A. Sestini are members of the INdAM Research group GNCS. The INdAM support through GNCS and Finanziamenti Premiali SUNRISE is gratefully acknowledged. A. Falini was also supported by the INdAM-GNCS project ``Finanziamenti Giovani Ricercatori 2019-2020''.
\end{acknowledgement}

\bibliographystyle{plain}
\bibliography{paper_3DBEM}

\end{document}